\documentclass[12pt,a4paper]{article}
\usepackage{amsmath}
\usepackage{amsfonts}
\usepackage{amsthm}
\usepackage{amscd}
\usepackage{amsbsy}
\usepackage{amsrefs}
\usepackage[latin2]{inputenc}
\usepackage{t1enc}
\usepackage[mathscr]{eucal}
\usepackage{indentfirst}
\usepackage{graphicx}
\usepackage{graphics}
\usepackage{pict2e}
\usepackage{epic}
\numberwithin{equation}{section}
\usepackage[margin=2.9cm]{geometry}
\usepackage{epstopdf}
\usepackage{tikz-cd}
\usepackage{chngcntr}
\usepackage{amssymb}
\usepackage{esint}
\usepackage{hyperref}
\usepackage{mwe}
\usepackage{mathrsfs}
\usepackage{xfrac}
\usepackage{verbatim}
\usepackage{mathpazo}

\theoremstyle{plain}
\newtheorem{Th}{Theorem}[section]

\newtheorem{Prop}[Th]{Proposition}

 \theoremstyle{definition}
\newtheorem{Def}[Th]{Definition}

\newtheorem{Rem}[Th]{Remark}
\newtheorem{?}[Th]{Problem}

\newcommand*\R{\mathbb{R}}

\newcommand*\Om{\Omega}

\newcommand{\Div}{\text{div}_x}
\newcommand{\Divh}{\text{div}_{x_h}}
\newcommand{\vectoru}{\mathbf{u}}

\newcommand{\vectorv}{\mathbf{v}}
\newcommand{\vectorV}{\mathbf{V}}
\newcommand{\vectorm}{\mathbf{m}}
\newcommand{\Epsilon}{\mathcal{E}}
\newcommand{\dx}{\text{ d}x}
\newcommand{\dt}{\text{ d}t\;}

\newcommand{\vectorphi}{\pmb{\varphi}}

\newcommand{\curlh}{\text{curl}_{x_h}}
\begin{document}

\title{Singular limits for compressible inviscid rotating fluids}

\author{Nilasis Chaudhuri
	\thanks{E-mail:\tt chaudhuri@math.tu-berlin.de}
}
\maketitle

\centerline{Technische Universit{\"a}t, Berlin}

\centerline{Institute f{\"u}r Mathematik, Stra\ss e des 17. Juni 136, D -- 10623 Berlin, Germany.}

\begin{abstract}
{We study singular limit for scaled barotropic Euler system modelling a rotating, compressible and inviscid fluid, where Mach and Rossby numbers are proportional to a small parameter $\epsilon$. If the fluid is confined to an infinite slab, the limit behaviour is identified as a horizontal motion of an incompressible inviscid system that is analogous to the Euler system. We consider a very general class of solutions, named \emph{dissipative solution} for the scaled compressible Euler systems and will show that it converges to a \emph{strong solution} of that incompressible inviscid system.}
\end{abstract}



\section{Introduction}
Several types of singular limits of Navier--Stokes system and Euler system have been studied in the last few decades. Here we devote ourselves on the models arising by the effect of rotation on fluids as described in
{Chemin et.al.} \cite{CDGG2006}.

Let $T>0$ and $\Omega \subset (\mathbb{R}^3) = \mathbb{R}^2 \times (0,1)$
be an infinite slab. We consider the scaled \emph{compressible Euler equation} in time-space cylinder $Q_T=(0,T)\times \Omega$ describing the time evolution of the mass density $\varrho=\varrho(t,x)$ and the momentum field $\vectorm=\vectorm(t,x)$ of a compressible rotating inviscid fluid with axis of rotation $\mathbf{b}=(0,0,1)$:

\begin{itemize}
	\item \textbf{Conservation of Mass: }
	\begin{align}
	\partial_t \varrho + \text{div}_x \vectorm&=0. \label{cee:cont}
	\end{align}
	\item \textbf{Conservation of Momentum:}
	\begin{align}
	\partial_t \mathbf{m} + \Div (\frac{ \mathbf{m} \otimes \mathbf{m}}{\varrho})+\frac{1}{\text{Ma}^2}\nabla_x p(\varrho)+\frac{1}{\text{Ro}} \mathbf{b} \times \mathbf{m}&=0. \label{cee:mom}
	\end{align}	
\end{itemize}

\begin{itemize}
	\item The scaled system contains charecteristic numbers:
	\subitem Ma-- Mach number,
	\subitem Ro-- Rossby number.\\
	Here  we consider,
	\begin{align}\label{scaling}
	\text{Ma} \approx \epsilon,\; \text{Ro} \approx \epsilon \text{ for } \epsilon>0.
	\end{align}

	\item \textbf{Pressure Law:} In an isentropic setting, the pressure $p$ and the density $\varrho$ of the fluid are interrelated by
		\begin{align}\label{p-condition}
		\begin{split}
		p(\varrho)=a \varrho^\gamma,\; a>0,\; \gamma> 1.
		\end{split}
		\end{align}
	
 	\item \textbf{Boundary condition:} Here we consider slip condition for velocity on the horizontal boundary i.e.
 	\begin{align}\label{BC}
 	\vectoru \cdot \mathbf{n}=0,\; \mathbf{n}=(0,0,\pm 1).
 	\end{align}
 	
 	\item \textbf{Far field condition:} Introducing the notation $x=(x_h,x_3)$ and $P_h(x)=x_h$ we assume this condition as,
 	\begin{align}\label{far_field}
 	\varrho \rightarrow \bar{\varrho}, \; \vectoru \rightarrow 0 \text{ as } \vert x_h \vert \rightarrow \infty.
 	\end{align}
 	\item \textbf{Initial data:}
 	For each $\epsilon>0$, we supplement the initial data as
 	\begin{align}\label{initial condition}
 	{\varrho_\epsilon(0,\cdot)=\varrho_{\epsilon, 0},\; \vectorm_\epsilon(0,\cdot)= \vectorm_{\epsilon, 0}.}
 	\end{align}
 	
 \end{itemize}

Our goal here is to study the effect of \emph{low Mach number limit} (also called \emph{incompressible limit}) and \emph{low Rossby number limit} on the system \eqref{cee:cont}-\eqref{cee:mom}. In low Mach number limit fluid becomes \emph{incompressible} and in low Rossby number limit indicates fast rotation of fluids  as a consequence of that fluid becomes \emph{planner} (two-dimensional).
 In our case when $\epsilon \rightarrow 0$, we have both phenomenon described above simultaneously, i.e the system \eqref{cee:cont}-\eqref{cee:mom} (\emph{primitive system}) which describes a compressible, rotating (3D) fluid is expected to become a system (\emph{target system}) that describes incompressible, planner(2D) fluid. \\

 There are mainly two approaches to deal with the singular limit problem.
 \begin{itemize}
 	\item[I.]  The first approach consists of considering a \emph{classical(strong)} solution of the \emph{primitive system} and expected that it converges to the classical solutions of the \emph{target system}. Here, the main and highly nontrivial issue is to ensure that the lifespan of the strong solutions is bounded below away from zero uniformly with
 	respect to the singular parameter.
 	
 	\item[II.] The second approach is based on the concept of \emph{weak} or \emph{measure--valued}  \emph{dissipative} solutions  of the \emph{primitive system}. Under proper choice of initial data one can show convergence 
	provided the \emph{target system} admits smooth solution.
 \end{itemize}
For the first approach in the low Mach number limits we have results by Schochet\cite{S1986}, Ebin \cite{E1977}, Kleinermann-Majda\cite{KM1981} and many others. For rotating fluids there are results by Babin et. al \cite{BMN1999}, \cite{BMN2001} and Chemin et. al. \cite{CDGG2006}. Since it is well known fact that strong solution of compressible inviscid fluids develop {singularities, the main difficulty here} is to show the existence time of \emph{primitive system} is independent of $\epsilon$.\\

In the case of second approach, a similar problem has been studied by Feireisl et. al. in \cite{FN2014}. They consider viscous fluid(Navier-Stokes System) and also high Reynolds number limit along with low Mach and low Rossby number limit. They have proved that weak solution of their system converges to the strong solution of the same \emph{target system} as we have . Relative energy inequality plays very important role in the proof. Also Feireisl et. al. \cite{FGN2012} have studied the low Mach and Rossby number limit for scaled Navier-Stokes system. In this case they obtain they obtain a similar system with some effect of viscous term.\\

The advantage(s) to consider the second approach is(are),
\begin{itemize}
	\item \emph{Weak} or \emph{measure valued} solutions to the \emph{primitive system} exist globally in time. Hence the result depends only on the life span of the \emph{target problem} that may be finite.
	\item {The result convergence holds for a large class of generalized solutions which indicates certain stability 
	of the limit solution of the target system.}
\end{itemize}

 Recently, the concept of measure--valued solutions has been studied in variuos context, like, analysis of numerical schemes etc. In the following articles \cite{AB1997}, \cite{PAW2015}, \cite{BF2018b}, \cite{B2018}, \cite{Bd2019}, \cite{FL2018} we observe the development of theory on measure valued solution for different models describing compressible fluids mainly with the help of Young measures. But Feireisl et. al. in \cite{FLM2019} and \cite{BeFH2019} give a new definition termed as \emph{dissipative solution} without involving Young measures. We will follow the last approach. \\

In  Feireisl et. al. \cite{FKM2019}, Bruell et. al. \cite{BrF2018} and Br{\'e}zina et. al.\cite{BM2018} show that \emph{measure valued solution} of \emph{primitive system} which describes some compressible inviscid fluid converges to \emph{strong solution} of incompressible \emph{target system}. \\

In second approach, it is very important to consider proper initial data mainly termed as \emph{well-prepared} and \emph{ill-prepared} initial data. Feireisl and Novotn{\'y} in \cite{FN2009b},  explain that for \emph{ill-prepared data} the presence of Rossby-acoustic waves play an important role in analysis of singular limits. Meanwhile this effect was absent in \emph{well-prepared} data. Here we deal with both types of initial data. \\

Our main goal is to prove that a \emph{dissipative solution} behaves similar to the weak solutions, i.e. in low Mach and low Rossby regime they converges to strong solution which descibes planner, incompressible fluid in 2D. Hence our plan for the article is,
\begin{itemize}
	\item[1.] Definition of dissipative solution.
	\item[2.] Singular limit for `\emph{well-prepared}' data.
	\item[3.] Singular limit for `\emph{ill-prepared}' data.
\end{itemize}

\subsection{Notation:}
\begin{itemize}
	\item To begin, we introduce a function
	$\chi = \chi(\varrho)$ such that
	\begin{align}
	\begin{split}
	\chi(\varrho) \in C_c^{\infty}	(0, \infty),\; 0 \leq \chi \leq 1,\; \chi(\varrho) = 1 \ \mbox{if}\ \frac{\bar{\varrho}}{2} \leq \varrho \leq 2 \bar{\varrho},\; \bar{\varrho}>0.
	\end{split}
	\end{align}
	For a function, $H = H(\varrho, \vectorm)$ we set
	\begin{align}
	\begin{split}
	[H]_{\text{ess}}= \chi(\varrho) H(\varrho, \vectorm),\; [H]_{\text{res}}= (1-\chi(\varrho)) H(\varrho, \vectorm).
	\end{split}
	\end{align}
	\item Following the article by Feireisl and Novotn{\'y} \cite{FN2014} we can conclude that the consideration of the domain $\R^2\times (0,1)$  is equivalent to consider $ \R^2 \times \mathbb{T}^1$. Since, the complete slip boundary conditions can be transformed to periodic ones by considering the space of symmetric functions, the horizontal components of the velocity are even and the vertical one odd with respect to the vertical variable, i.e.
	\begin{align}\label{Symmetry_class}
	\begin{split}
	\varrho(t,x_h,-x_3)= \varrho(t,x_h,x_3),\; &\vectorv_h(t,x_h,-x_3)=\vectorv_h(t,x_h,x_3),\\
	&\; v_3(t,x_h,-x_3)=-v_3(t,x_h,x_3).
	\end{split}
	\end{align}
	for all $t\in (0,T),\; x_h \in \R^2,\; x_3 \in \mathbb{T}^{1}$. This equivalence has been described in Ebin\cite{E1983}. So from now on we consider $\Om= \R^2\times \mathbb{T}^1$.
	\item Let us define pressure potential as,
	\begin{align}\label{Pressure_potential}
	P(\varrho)=\varrho \int_{\bar{\varrho}}^{\varrho} \frac{p(z)}{z^2} \; \text{d}z.
	\end{align}
	As a consequence of that we have,
	\begin{align}
	\varrho P^{\prime}(\varrho) -P(\varrho)= p(\varrho) \text{ and } \varrho P^{\prime \prime} (\varrho)=p^{\prime}(\varrho) \text{ for } \varrho>0.
	\end{align}
	
\end{itemize}

\section{Definition of dissipative solution}

We first want to give a definition of \emph{dissipative solution} for our choosen system. Basari{\'c} in \cite{Bd2019} introduces a dissipative measure valued solution for Euler equation with damping by introducing Young measures in an unbounded domain. For bounded domain,  Feireisl et. al.\cite{FLM2019} and \cite{BeFH2019} recently introduced a concept of dissipative solution without defining solution through Young Measures. We follow the second approach to define \emph{dissipative solution}.

\begin{Def}\label{dmv_defn}
	Let $\epsilon>0$ and $\bar{\varrho}>0$. We say functions $\varrho_{\epsilon},\vectorm_{\epsilon}$  with,
	\begin{align}
	\varrho_{\epsilon}-\bar{\varrho} \in  C_{\text{weak}}([0,T];L^2+L^{\gamma}(\Om)),\;\varrho_\epsilon \geq 0,\; \vectorm_{\epsilon} \in   C_{\text{weak}}([0,T];L^2+L^{\frac{2\gamma}{\gamma+1}}(\Om))
	\end{align}
	and also satisfying \eqref{Symmetry_class} are a \emph{dissipative solution} to the compressible Euler equation \eqref{cee:cont}-\eqref{far_field} with initial data $\varrho_{0,\epsilon}, \vectorm_{0,\epsilon}$ satisfying,
	\begin{align}\label{fe_id}
	\begin{split}
	\varrho_{0,{\epsilon}}\geq 0,\; E_{0,{\epsilon}}=\int_{ \Om} \bigg( \frac{1}{2} \frac{\vert \vectorm_{0,{\epsilon}} \vert^2}{\varrho_{0,{\epsilon}}} + P(\varrho_{0,{\epsilon}})-(\varrho_{0,{\epsilon}}-\bar{\varrho})P^{\prime}(\bar{\varrho})\bigg)\dx < \infty,
	\end{split}
	\end{align}
	if there exist the \emph{turbulent defect measures}
	\begin{align}
	\begin{split}
	&{\mathfrak{R}_{v_{\epsilon}} \in L^{\infty}(0,T;\mathcal{M}^+({\Om;\mathbb{R}^{d\times d}_{\text{sym}}}))},\; \mathfrak{R}_{p_{\epsilon}} \in L^{\infty}(0,T;\mathcal{M}^{+}({\Om})),
	\end{split}
	\end{align}
	such that the following holds,
	\begin{itemize}
		\item \textbf{Equation of Continuity:}
		For any $\tau \in (0,T) $ and any $\varphi \in C_{c}^{1}([0,T]\times \bar{\Om})$ {it holds}
		\begin{align}\label{continuity_eqn_weak}
		\begin{split}
		&\big[ \int_{ \Om}{( \varrho_{\epsilon}-\bar{\varrho})} \varphi \dx\big]_{t=0}^{t=\tau}=
		\int_0^{\tau} \int_{\Om} [ \varrho_{\epsilon} \partial_t \varphi +  \vectorm_{\epsilon} \cdot \nabla_x \varphi] \dx \dt .
		\end{split}
		\end{align}
		
		\item \textbf{Momentum equation:} For any $\tau\in (0,T)$ and any $\pmb{\varphi} \in C^{1}_c([0,T]\times \Om;\mathbb{R}^d)$, it holds
		\begin{align}\label{momentum_eqn_weak}
		\begin{split}
		&\bigg[\int_{\Om}  \vectorm_{\epsilon}(\tau,\cdot)\cdot \vectorphi(\tau,\cdot) \dx\bigg]_{t=0}^{t=\tau} \\
		&=\int_0^{\tau}\int_{\Om} [ \vectorm_{\epsilon}\cdot \partial_{t} \vectorphi + (\frac{ \vectorm_{\epsilon} \otimes \vectorm_{\epsilon}}{\varrho_{\epsilon}}) : \nabla_x \vectorphi + p(\varrho_{\epsilon}) \Div \vectorphi + \mathbf{b} \times \mathbf{m}_{\epsilon} \cdot \vectorphi ] \dx \dt \\
		&+ \int_0^{\tau}\int_{\Om} \nabla_x \vectorphi : \text{d}[\mathfrak{R}_{p_{\epsilon}} +\mathfrak{R}_{p_{\epsilon}} \mathbb{I}]  \dt .
		\end{split}
		\end{align}
		\item \textbf{Energy inequality:}  The total energy $E$ is defined in $[0,T)$ as,
		\begin{align*}
		E_{\epsilon}(\tau)= \int_{\Om} &\bigg( \frac{1}{2} \frac{\vert \vectorm_{\epsilon} \vert^2}{\varrho_{\epsilon}} +(P(\varrho_{\epsilon})-(\varrho_{\epsilon}-\bar{\varrho})P^{\prime}(\bar{\varrho})) \bigg)(\tau ,\cdot)\dx.
		\end{align*}
			It satisfies,
		\begin{align}\label{energy_inequality}
		\begin{split}
		\bigg[ \psi \bigg( E_{\epsilon}(\tau) + \int_{ \Om}\bigg[ \frac{1}{2} \text{d}\text{ trace} \mathfrak{R}_{v_{\epsilon}}(\tau,\cdot)+ \frac{1}{\gamma-1} \text{d}\mathfrak{R}_{p_{\epsilon}}(\tau,\cdot) \bigg] \bigg) \bigg]_{\tau=\tau_1-}^{\tau=\tau_2+} \\ \leq \int_{\tau_1}^{\tau_2} \partial_t(\psi) \bigg( E_{\epsilon}(t)+ \int_{ \bar{\Om}} \frac{1}{2} \text{d} \text{ trace} \mathfrak{R}_{{v_\epsilon}}+\int_{ \bar{\Om}} \frac{1}{\gamma-1} \text{d}  \mathfrak{R}_{{p_\epsilon}}\bigg) \dt
		\end{split}
		\end{align}
		for any  $0\leq \tau_1\leq \tau_2<T $ and any $\psi \in C_{c}^{1}[0,T)$, $\psi\geq 0$.		
	\end{itemize}
\end{Def}

\begin{Rem}
 In \eqref{energy_inequality}, the initial energy satisfies 	$E_{\epsilon}(0-)=E_{0,\epsilon}$.
\end{Rem}

\begin{Th}
	Suppose $\Om$ be the domain specified above and pressure follows \eqref{p-condition}. If $(\varrho_{0,\epsilon},\vectorm_{0,\epsilon})$ satisfies \eqref{fe_id}, then there exists \emph{dissipative solution} as defined in definition \eqref{dmv_defn}.
\end{Th}

The proof this theorem follows in similar lines of Breit, Feireisl and Hofmanov{\'a} as in \cite{BeFH2019} and \cite{BeFH2019(2)}. We have to adopt it for unbounded domain as suggested in Basari{\'c}\cite{Bd2019}.\\

\section{Identification of \emph{Target system}}

Taking motivation from \cite{FN2014} we expect the target system as,
\begin{align}\label{limit_system}
\begin{split}
&\frac{p^{\prime}(\bar{\varrho})}{\bar{\varrho}}\nabla_x q +  \mathbf{b} \times \vectorv =0,\\
&\partial_{t}(\Delta_{x_h}q-\frac{1}{p^{\prime}(\bar{\varrho})}q) + (\nabla_{x_h}^{\perp}q\cdot \nabla_{x_h})(\Delta_{x_h}q- \frac{1}{p^{\prime}(\bar{\varrho})}q)=0, \; \text{ in }\mathbb{R}^2.
\end{split}
\end{align}

\subsection{Informal justification:}
Here is an informal justification how we obtain the \emph{target system} as in \eqref{limit_system}.
First we note that $(\bar{\varrho},0)$ is a \emph{steady state} solution for \eqref{cee:cont}-\eqref{far_field}.

Let us consider
\begin{align*}
&\varrho_{\epsilon}= \bar{\varrho} + \epsilon \varrho_{\epsilon}^{(1)}+ \epsilon^2  \varrho_{\epsilon}^{(2)}+ \cdots,\\
&\vectorm_\epsilon= \bar{\varrho}\vectorv + \epsilon  \vectorm_\epsilon^{(1)} + \epsilon^2 \vectorm_\epsilon^{(2)} + \cdots.
\end{align*}
As a consequence of the above we obtain,
\begin{align*}
p(\varrho_{\epsilon})= p(\bar{\varrho}) + \epsilon p^{\prime}(\bar{\varrho}) \varrho_{\epsilon}^{(1)} + \epsilon^2 (  p^{\prime}(\bar{\varrho})\varrho_{\epsilon}^{(2)}+  \frac{1}{2}p^{\prime\prime}(\bar{\varrho})(\varrho_{\epsilon}^{(1)})^{2}) + o(\epsilon^3).
\end{align*}

So we obtain,
\begin{align*}
\bar{\varrho}\Div \vectorv + \epsilon (\partial_{t} \varrho_{\epsilon}^{(1)}+ \Div \vectorm_\epsilon^{(1)}) + o(\epsilon^2)=0
\end{align*}
and
\begin{align*}
&\bar{\varrho}(\partial_{t} \vectorv + (\vectorv\cdot \nabla_x)\vectorv) + \nabla_x (  p^{\prime}(\bar{\varrho})\varrho_{\epsilon}^{(2)}+  \frac{1}{2}p^{\prime\prime}(\bar{\varrho})(\varrho_{\epsilon}^{(1)})^{2}) + \mathbf{b} \times \vectorm_\epsilon^{(1)}\\
&+ \frac{1}{\epsilon} (p^{\prime}(\bar{\varrho}) \nabla_x \varrho_{\epsilon}^{(1)}+ \bar{\varrho}\mathbf{b} \times \vectorv) + o(\epsilon)=0.
\end{align*}

Further we assume $ \big( \frac{\varrho_\epsilon - \bar{\varrho}}{\epsilon} \big) \rightarrow q$ and $ \vectorm_\epsilon\rightarrow \bar{\varrho} \vectorv$ in some strong sense. Then as a consequence we have
\begin{align*}
&p^{\prime}(\bar{\varrho})\nabla_x q +  \bar{\varrho}\mathbf{b} \times \vectorv =0,\\
&\Div \vectorv =0 ,\\
&\bar{\varrho}(\partial_{t} \vectorv + (\vectorv \cdot \nabla_{x})\vectorv) + \nabla_{x} \Pi_1+ \mathbf{b} \times \vectorm_1=0,\\
& \partial_{t} q + \Div \vectorm_1=0.
\end{align*}

Few consequence from above relations,
\begin{align*}
q_{x_3}=0,\; q(x)=q(x_h), \nabla_{x_h}^{\perp}q=\frac{\bar{\varrho}}{p^{\prime}(\bar{\varrho})}\vectorv_h, \vectorv_h=(v_1,v_2),\\
\Divh \vectorv_h=0, v_{3_{x_3}}=0, \text{ here } \nabla_{x_h}^{\perp}\equiv (-\frac{\partial}{\partial x_2},\frac{\partial}{\partial x_1}).
\end{align*}
Assuming smoothness,
\begin{align*}
v_{1_{x_3}}=0,\; v_{2_{x_3}}=0
\end{align*}
So,
\begin{align*}
\vectorv(x)=\vectorv(x_h).
\end{align*}
Also boundary condition $\frac{\vectorm_\epsilon}{\varrho_{\epsilon}}\cdot \mathbf{n} =0 $ will lead to conclude
\begin{align*}
v_3(x_h,x_3)=0.
\end{align*}
Then we will have $ \vectorv=(\vectorv_h(x_h),0)$. Finally we obtain,
\begin{align}
\begin{split}
&\frac{p^{\prime}(\bar{\varrho})}{\bar{\varrho}}\nabla_{x_h} q +  \mathbf{b} \times \vectorv =0,\\
&\partial_{t}(\Delta_{x_h}q-\frac{1}{p^{\prime}(\bar{\varrho})}q) + (\nabla_{x_h}^{\perp}q\cdot \nabla_{x_h})(\Delta_{x_h}q- \frac{1}{p^{\prime}(\bar{\varrho})}q)=0.
\end{split}
\end{align}
Here $q$ can be viewed as a kind of \emph{stream function}. As a simple consequence of above we have, $\Divh \vectorv=0$.

\subsection{Reguarity of Target System \eqref{limit_system}}
Since the \eqref{limit_system} is possesses the same structure as 2D Euler equations. We expect solution to be as regular as the initial data and exists globally in time. In particular, we may use the abstract theory of Oliver\cite{O1997}, Theorem 3, to obtain the result:

\begin{Prop}\label{target_sys}
	Suppose that
	\begin{align*}
	q_0\in W^{m,2}(\mathbb{R}^2) \text{ for } m\geq 4.
	\end{align*}
	Then, the problem \eqref{limit_system}  admits a solution $q$, unique in the class
	\begin{align*}
	q \in C([0,T];W^{m,2}(\R^2)) \cap C^{1}([0,T];W^{m-1,2}(\R^2)).
	\end{align*}	
\end{Prop}

\section{ Singular limit for "Well-prepared" initial Data }
After giving a proper definition of \emph{well-prepared} data we will state our mail result and will prove it in the following sections.
\subsection{ Definition of "Well-prepared Data"}
\begin{Def}
	We say that the set of initial data $(\varrho_{0,\epsilon},\vectorm_{0,\epsilon})_{(\epsilon>0)} $ is "well-prepared" if,
	\begin{align}\label{well_prepared_id}
	\begin{split}
	&\varrho_{0,\epsilon}=\bar{\varrho}+\epsilon \varrho_{0,\epsilon}^{(1)},\; \{\varrho_{0,\epsilon}\}_{\{\epsilon>0\} }\text{ is bounded in } L^2\cap L^{\infty}(\Om), \; \varrho_{0,\epsilon}^{(1)}\rightarrow q_0 \text{ in } L^2(\Om),\\
	&\frac{\vectorm_{0,\epsilon}}{\varrho_{0,\epsilon}}\rightarrow  \vectorv_{0} = \big(\vectorv_{0}^{(1)},\vectorv_{0}^{(2)},0) \text{ in } L^2(\Omega;\R^3\big) \text{ with the following relation}\\
	& -\Delta_{x_h} q_0 = \bar{\varrho}  \curlh P_h(\vectorv_0) .
	\end{split}
	\end{align}
\end{Def}

For $\epsilon>0$, $(\varrho_\epsilon, \vectorm_\epsilon)$ are solutions of \eqref{cee:cont}-\eqref{far_field} with initial data \eqref{initial condition}. If we have "Well--prepared initial data" then claim is $ (\varrho_\epsilon, \vectorm_\epsilon) \rightarrow (\bar{\varrho},\vectorv_h,q)$ where $(\bar{\varrho},\vectorv_h,q) $ is a classical solution of the \eqref{limit_system} with initial data $(P_h(\vectorv_{0}),0)$ in $\R^2$.

\subsection{Main Theorem:}

\begin{Th}\label{th1}
	Let pressure $p$ follows \eqref{p-condition}. Further we assume that the initial data is \emph{well-prepared}, i.e. it follows \eqref{well_prepared_id}. Then after taking a subsequence, the following holds,
	\begin{align}
	\begin{split}
	&\varrho_{\epsilon}^{1} \equiv \frac{\varrho_{\epsilon}-\bar{\varrho}}{\epsilon} \rightarrow q \text{ weakly-(*) in } L^{\infty}(0,T;L^2+L^\gamma(\Om)),\\
	&\vectorm_\epsilon \rightarrow \bar{\varrho} \vectorv \text{ weakly-(*) in } L^{\infty}(0,T;L^2+L^{\frac{2\gamma}{\gamma+1}}(\Om)).\\
	\end{split}
	\end{align}
	Now if we assume that $\vectorv_0\in W^{k,2}$ with $k\geq 3$. Then $(q,\vectorv)$ solves \eqref{limit_system} with initial data $q(0,\cdot)=q_0$ where $q_0$ satisfies \eqref{well_prepared_id}.	
\end{Th}
In the following subsections we will complete the proof.

\subsection{Relative energy inequality}
Relative energy plays an important role to conclude the proof of the theorem \eqref{th1}.
For $\epsilon=1$, the \emph{relative energy} defines as,
\begin{align}\label{rel-ent}
\Epsilon(t)=\Epsilon(\varrho,\vectorm \vert \tilde{\varrho},\tilde{\vectoru})(t):= \int_{\Omega}\frac{1	}{2} \varrho \vert \frac{\vectorm}{\varrho}-\tilde{\vectoru}\vert^2 + (P(\varrho)-P(\tilde{\varrho}) -P^{\prime}(\tilde{\varrho})(\varrho -\tilde{\varrho})) (t,\cdot) \dx ,
\end{align}
where $\tilde{\varrho},\tilde{\vectoru}$ are arbitrary smooth test functions with $\tilde{\varrho}-\bar{\varrho}$ and $\tilde{\vectoru}$ have compact support .

\begin{Rem}
	The relative energy is a coercive functional (see. Bruell et. al. \cite{BrF2018}) satisfying the estimate,
	\begin{align}
	\begin{split}
	 \Epsilon(\varrho,\vectorm \vert \tilde{\varrho},\tilde{\vectoru})(t) \geq& \int_{ \Om} \bigg[  \vert \frac{\vectorm}{\varrho}-\tilde{\vectoru}\vert^2  \bigg]_{\text{ess}} \dx + \int_{ \Om} \bigg[  \frac{\vert \vectorm \vert^2}{\varrho}  \bigg]_{\text{res}} \dx \\
	\quad & + \int_{ \Om} [(\varrho-\tilde{\varrho})^2]_{\text{ess}}\dx + \int_{ \Om} [1]_{\text{res}} + [\varrho^{\gamma}]_{\text{res}}\dx.
	\end{split}
	\end{align}
\end{Rem}

From \eqref{rel-ent} we can deduce that,
\begin{align}
\begin{split}
\Epsilon(\tau) &=\int_{ \Om} \big( \frac{1}{2} \frac{\vert \vectorm \vert^2}{\varrho}+ P(\varrho) \big) \dx - \int_{ \Om}  \vectorm  \cdot \tilde{\vectoru} \dx \\
& \quad +\int_{ \Om} \frac{1}{2} \varrho  \vert \tilde{\vectoru} \vert^2 \dx - \int_{ \Om}\varrho P^\prime(\tilde{\varrho}) \dx + \int_{ \Om} p(\tilde{\varrho}) \dx \\
&= \int_{\Om} \bigg( \frac{1}{2} \frac{\vert \vectorm \vert^2}{\varrho} +P(\varrho)-(\varrho-\bar{\varrho})P^{\prime}(\bar{\varrho}) \bigg)\dx\\
& \quad - \int_{ \Om}  \vectorm  \cdot \tilde{\vectoru} \dx +\int_{ \Om} \frac{1}{2} \varrho  \vert \tilde{\vectoru} \vert^2 \dx - \int_{ \Om}\varrho P^\prime(\tilde{\varrho}) \dx + \int_{ \Om} p(\tilde{\varrho}) \dx \\
& \quad + \int_{ \Om} (\varrho-\bar{\varrho})P^{\prime}(\bar{\varrho}) \dx.
\end{split}
\end{align}
Now using the \emph{dissipative solution} we obtain,
\begin{align*}
\Epsilon(\tau)+&  \int_{ \Om}\bigg[ \frac{1}{2} \text{d}\text{ trace} \mathfrak{R}_{v}(\tau,\cdot)+ \frac{1}{\gamma-1} \text{d}\mathfrak{R}_{p}(\tau,\cdot) \bigg]\\
\quad &\leq E_0 - \int_{ \Om} \vectorm_0 \tilde{\varrho}(0,\cdot)\dx + \int_{ \Om} \frac{1}{2} \varrho_{0} \vert \tilde{\vectoru}(0,\cdot) \vert^2 \dx - \int_{ \Om}\varrho_{0}P^{\prime}(\tilde{\varrho}(0,\cdot)) \\
&\quad +\int_{ \Om} (\varrho_{0}-\bar{\varrho}) P^{\prime}(\bar{\varrho}) \dx + \int_{ \Om} p(\tilde{\varrho}(0,\cdot))\dx\\
\quad & \quad -\int_0^{\tau}\int_{\Om}  \vectorm\cdot \partial_{t} \tilde{\vectoru} - (\frac{ \vectorm \otimes \vectorm}{\varrho}) : \nabla_x \tilde{\vectoru} - p(\varrho) \Div \tilde{\vectoru} + \mathbf{b} \times \mathbf{m} \cdot \tilde{\vectoru}  \dx \dt \\
& \quad - \int_0^{\tau}\int_{\Om} \nabla_x \tilde{\vectoru} : \text{d}[\mathfrak{R}_{v} +\mathfrak{R}_{p} \mathbb{I}]  (t,\cdot)\dt \\
\quad & \quad + \int_0^{\tau} \int_{\Om}  \varrho \partial_t  (\frac{1}{2} \vert \tilde{\vectoru} \vert^2) +  \vectorm \cdot \nabla_x ( \frac{1}{2}   \vert \tilde{\vectoru} \vert^2) \dx \dt \\
&\quad -\int_0^{\tau} \int_{\Om}  \varrho \partial_t   P^\prime(\tilde{\varrho})  -  \vectorm \cdot \nabla_x  P^\prime(\tilde{\varrho}) \dx \dt + \int_0^{\tau} \int_{\Om} p^{\prime}(\tilde{\varrho}) \partial_{t}\tilde{\varrho} \dx \dt.
\end{align*}

Finally we obtain,

\begin{align*}
	\Epsilon(\tau)+&  \int_{ \Om}\bigg[ \frac{1}{2} \text{d}\text{ trace} \mathfrak{R}_{v}(\tau,\cdot)+ \frac{1}{\gamma-1} \text{d}\mathfrak{R}_{p}(\tau,\cdot) \bigg]\\
	\quad &\leq \Epsilon(0)  -\int_0^{\tau}\int_{\Om}  (\vectorm-\varrho \tilde{\vectoru})\cdot \partial_{t} \tilde{\vectoru}\dx \dt\\
	&\quad  - \int_0^{\tau}\int_{\Om}(\frac{ (\vectorm-\varrho \tilde{\vectoru}) \otimes \vectorm}{\varrho}) : \nabla_x \tilde{\vectoru} \dx \dt  \\
	&\quad  - \int_0^{\tau}\int_{\Om} (p(\varrho)-p(\tilde{\varrho})) \Div \tilde{\vectoru} \dx \dt \\
	& \quad +\int_0^{\tau}\int_{\Om} \mathbf{b} \times \mathbf{m} \cdot ( \tilde{\vectoru}- \frac{\vectorm}{\varrho})  \dx \dt +\int_0^{\tau}\int_{\Om} (\tilde{\varrho}-\varrho) \partial_{t} P^{\prime}(\tilde{\varrho}) \dx\dt \\
	&\quad + \int_0^{\tau}\int_{\Om}(\tilde{\varrho}\tilde{\vectoru}- \vectorm) \nabla_x P^{\prime}(\tilde{\varrho}) \dx \dt- \int_0^{\tau}\int_{\Om} \nabla_x \tilde{\vectoru} : \text{d}[\mathfrak{R}_{v} +\mathfrak{R}_{p} \mathbb{I}]  (t,\cdot)\dt . \\
\end{align*}

Now we consider the scaling on \eqref{rel-ent} and rephrase it as,
\begin{align}\label{rel-ent-epsilon}
\begin{split}
\Epsilon_{\epsilon}(t)&=\Epsilon(\varrho_{\epsilon},\vectorm_{\epsilon} \vert \tilde{\varrho},\tilde{\vectoru})(t)\\
&:= \int_{\Omega}\big[ \frac{1	}{2} \varrho_{\epsilon} \vert \frac{\vectorm}{\varrho_{\epsilon}}-\tilde{\vectoru}\vert^2 +\frac{1}{\epsilon^2} (P(\varrho_{\epsilon})-P(\tilde{\varrho}) -P^{\prime}(\tilde{\varrho})(\varrho_{\epsilon} -\tilde{\varrho})) \big] (t,\cdot) \dx .
\end{split}
\end{align}
We write the scaled version of the inequality as,
\begin{align}\label{rel-ent-m}
\begin{split}
	\Epsilon_{\epsilon}(\tau)&+  \int_{ \Om}\bigg[ \frac{1}{2} \text{d}\text{ trace} \mathfrak{R}_{v_{\epsilon}}(\tau,\cdot)+ \frac{1}{\gamma-1} \text{d}\mathfrak{R}_{p_{\epsilon}}(\tau,\cdot) \bigg]\\
	\quad &\leq \Epsilon_{\epsilon}(0)  -\int_0^{\tau}\int_{\Om}  (\vectorm_{\epsilon}-\varrho_{\epsilon} \tilde{\vectoru})\cdot \partial_{t} \tilde{\vectoru}\dx \dt\\
	&\quad  - \int_0^{\tau}\int_{\Om}(\frac{ (\vectorm_{\epsilon}-\varrho_{\epsilon} \tilde{\vectoru}) \otimes \vectorm_{\epsilon}}{\varrho_{\epsilon}}) : \nabla_x \tilde{\vectoru} \dx \dt  \\
	&\quad  - \frac{1}{\epsilon^2} \int_0^{\tau}\int_{\Om} (p(\varrho_{\epsilon})-p(\tilde{\varrho})) \Div \tilde{\vectoru} \dx \dt \\
	& \quad +\frac{1}{\epsilon} \int_0^{\tau}\int_{\Om} \mathbf{b} \times \mathbf{m}_{\epsilon} \cdot ( \tilde{\vectoru}- \frac{\vectorm_{\epsilon}}{\varrho_{\epsilon}})  \dx \dt  +\frac{1}{\epsilon^2} \int_0^{\tau}\int_{\Om} (\tilde{\varrho}-\varrho_{\epsilon}) \partial_{t} P^{\prime}(\tilde{\varrho}) \dx\dt \\
	&\quad + \frac{1}{\epsilon^2} \int_0^{\tau}\int_{\Om}(\tilde{\varrho}\tilde{\vectoru}- \vectorm_{\epsilon}) \nabla_x P^{\prime}(\tilde{\varrho}) \dx \dt - \int_0^{\tau}\int_{\Om} \nabla_x \tilde{\vectoru} : \text{d}[\mathfrak{R}_{v} +\mathfrak{R}_{p} \mathbb{I}]  (t,\cdot)\dt . \\
\end{split}
\end{align}
We rewrite the above inequality as,
\begin{align}
\begin{split}
\big[\Epsilon_{\epsilon}(\varrho_{\epsilon},\vectorm_{\epsilon}| \tilde{\varrho};\tilde{\vectoru})\big]_{0}^{\tau } &+ \int_{ \Om}\bigg[ \frac{1}{2} \text{d}\text{ trace} \mathfrak{R}_{v_{\epsilon}}(\tau,\cdot)+ \frac{1}{\gamma-1} \text{d}\mathfrak{R}_{p_{\epsilon}}(\tau,\cdot) \bigg]\\
&+ \big[\mathcal{R}_{\epsilon}(\varrho_{\epsilon},\vectorm_{\epsilon}| \tilde{\varrho};\tilde{\vectoru})\big] \leq 0.
\end{split}
\end{align}

Suppose, $(\tilde{\varrho}-\bar{\varrho}, \tilde{\vectoru})\in C^1([0,T]; H^m(\Om))\times C^1([0,T]; H^m(\Om))$  with $m$ large. We know $C^1([0,T];C_{c}^{\infty}(\Om))$ is dense in Sobolev space $C^1([0,T]; H^m(\Om))$. For $\delta>0$ we have, $\tilde{r} \in C^1([0,T];C_{c}^{\infty}(\Om))$ and $\tilde{\vectorv} \in C^1([0,T];C_{c}^{\infty}(\Om))$.
\begin{align*}
\Vert \tilde{r} - \tilde{\varrho} \Vert_{ C^1([0,T]; H^m(\Om))} + \Vert \tilde{\vectorv} - \tilde{\vectoru} \Vert_{ C^1([0,T]; H^m(\Om))} < \delta.
\end{align*}
Following Theorem 2.3 of \cite{Bd2019} we can show that,
\begin{align}
\begin{split}
& \big[\Epsilon_{\epsilon}(\varrho_{\epsilon},\vectorm_{\epsilon}| \tilde{\varrho};\tilde{\vectoru})\big]_{0}^{\tau } + \int_{ \Om}\bigg[ \frac{1}{2} \text{d}\text{ trace} \mathfrak{R}_{v_{\epsilon}}(\tau,\cdot)+ \frac{1}{\gamma-1} \text{d}\mathfrak{R}_{p_{\epsilon}}(\tau,\cdot) \bigg]+ \big[\mathcal{R}_{\epsilon}(\varrho_{\epsilon},\vectorm_{\epsilon}| \tilde{\varrho};\tilde{\vectoru})\big]\\
&\leq \big[\Epsilon_{\epsilon}(\varrho_{\epsilon},\vectorm_{\epsilon}| \tilde{r};\tilde{\vectorv})\big]_{0}^{\tau } + \int_{ \Om}\bigg[ \frac{1}{2} \text{d}\text{ trace} \mathfrak{R}_{v_{\epsilon}}(\tau,\cdot)+ \frac{1}{\gamma-1} \text{d}\mathfrak{R}_{p_{\epsilon}}(\tau,\cdot) \bigg]\\
&\quad +\big[\mathcal{R}_{\epsilon}(\varrho_{\epsilon},\vectorm_{\epsilon}| \tilde{r};\tilde{\vectorv})\big] + C \delta \\
& \leq C\delta.\\
\end{split}
\end{align}
Thus for Sobolev functions, relative energy inequality \eqref{rel-ent-m} is true.
\subsection{Convergence: Part 1}

First with $\tilde{\vectoru}= 0,\; \tilde{\varrho}=\bar{\varrho}$ as test functions we obtain the following bounds,
\begin{align*}
&\text{ess}\sup_{t\in (0,T)} \Vert \frac{\vectorm_{\epsilon}}{\sqrt{\varrho_{\epsilon}}} \Vert_{L^2(\Om;\R^3)} \leq C,\\
&\text{ess}\sup_{t\in (0,T)} \bigg\Vert \bigg[\frac{\varrho_{\epsilon}-\bar{\varrho}}{\epsilon}\bigg]_{\text{ess}}\bigg\Vert_{L^2(\Om)} \leq C,\\
&\text{ess}\sup_{t\in (0,T)} \Vert [\varrho_{\epsilon}]_{\text{res}} \Vert_{L^\gamma(\Om)}^{\gamma}  + \text{ess}\sup_{t\in (0,T)} \Vert [1]_{\text{res}} \Vert_{L^\gamma(\Om)}^{\gamma} \leq \epsilon^2 C .
\end{align*}

As an immediate consequence of above we have,
\begin{align*}
&\varrho_{\epsilon}^{(1)} \equiv \bigg(\frac{\varrho_{\epsilon}-\bar{\varrho}}{\epsilon}\bigg) \rightarrow \varrho^{(1)} \text{ wealy-}(*) \text{ in } L^{\infty}(0,T;L^2+L^{\gamma}(\Om)),\\
&\frac{\vectorm_{\epsilon}}{\sqrt{\varrho_{\epsilon}}} \rightarrow \hat{\vectorm} \text{ wealy-}(*) \text{ in } L^{\infty}(0,T;L^2(\Om;\R^3)).
\end{align*}
We obtain,
\begin{align*}
\varrho_{\epsilon} \rightarrow \bar{\varrho} \text{ in } L^{\infty}(0,T;L^2+L^{\gamma}(\Om))
\end{align*}
and
\begin{align*}
\Vert \vectorm_{\epsilon} \Vert_{L^{\infty}(0,T;L^2+L^{\sfrac{2\gamma}{\gamma+1}}(\Om;\R^3))} \leq C.
\end{align*}
Thus we have,
\begin{align*}
\vectorm_{\epsilon}\rightarrow \vectorm \text{ wealy-}(*) \text{ in } {L^{\infty}(0,T;L^2+L^{\sfrac{2\gamma}{\gamma+1}}(\Om;\R^3))}.
\end{align*}

Now from these two we write,
$\vectorm= \bar{\varrho} \vectoru$ and
letting $\epsilon \rightarrow 0$ in the continuity equation we have,
\begin{align*}
\bar{\varrho}\int_0^{\tau} \int_{\Om}  \vectoru \cdot \nabla_x \varphi \dx \dt=0.
\end{align*}
Furthermore multiplying momentum equation by $ \epsilon$ and letting $\epsilon \rightarrow 0$ we get the diagonasticequation,
\begin{align}
\mathbf{b} \times \vectoru + \frac{p^{\prime}(\bar{\varrho})}{\bar{\varrho}} \nabla_{x} \varrho^{1}=0,
\end{align}
in the sense of distributions.\\
Clearly from last relation we have $\varrho^{(1)}$ is independent of $x_3$, i.e. $\varrho^{(1)}=\varrho^{(1)}(x_{h})$.
Thus,
\begin{align*}
\vectoru=\vectoru(x_h),\; \Div \vectoru = \Divh \vectoru_{h}=0.
\end{align*}
Now from the symmetry class in which $\vectoru$ belongs we have,
\begin{align}
u_3=0, \; \vectoru=[\vectoru_h,0].
\end{align}
\subsection{Convergence: Part 2}

We recall the target system here,

\begin{align}
\begin{split}
&\frac{p^{\prime}(\bar{\varrho})}{\bar{\varrho}}\nabla_x q +  \mathbf{b} \times \vectorv =0,\\
&\partial_{t}(\Delta_{x_h}q-\frac{1}{p^{\prime}(\bar{\varrho})}q) + (\nabla_{x_h}^{\perp}q\cdot \nabla_{x_h})(\Delta_{x_h}q- \frac{1}{p^{\prime}(\bar{\varrho})}q)=0.
\end{split}
\end{align}
Let $(q,\vectorv)$ be solution of the above system with initial data $q_0,\vectorv_0$ satisfying \eqref{well_prepared_id}. Our goal is to show that $(\varrho^{(1)},\vectoru) \equiv  (q,\vectorv)$.
Here we choose proper test functions and will show that $\lim\limits_{\epsilon\rightarrow 0} \Epsilon_{\epsilon}(t)=0$.\\
We consider,
\begin{align}
\tilde{\vectoru}= \vectorv,\; \tilde{\varrho}=\bar{\varrho}+ \epsilon q.
\end{align}

We rewrite the relative energy inquality as,
 \begin{align*}
 \Epsilon_{\epsilon}(\tau)+&  \int_{ \Om}\bigg[ \frac{1}{2} \text{d}\text{ trace} \mathfrak{R}_{v_{\epsilon}}(\tau,\cdot)+ \frac{1}{\gamma-1} \text{d}\mathfrak{R}_{p_{\epsilon}}(\tau,\cdot) \bigg] \\
 \quad &\leq \Epsilon_{\epsilon}(0)  -\int_0^{\tau}\int_{\Om}  (\vectorm_{\epsilon}-\varrho_{\epsilon} \vectorv )\cdot ( \partial_{t} \vectorv  + (\vectorv\cdot \nabla_{x_h})\vectorv )\dx \dt\\
 &\quad  - \int_0^{\tau}\int_{\Om}(\frac{ (\vectorm_{\epsilon}-\varrho_{\epsilon} \vectorv) \otimes( \vectorm_{\epsilon}-\varrho_{\epsilon}\vectorv)}{\varrho_{\epsilon}}) : \nabla_{x}\vectorv  \dx \dt  \\
 &\quad  - \frac{1}{\epsilon^2} \int_0^{\tau}\int_{\Om} (p(\varrho_{\epsilon})-p(\tilde{\varrho})) \Div  \vectorv \dx \dt \\
 & \quad +\frac{1}{\epsilon} \int_0^{\tau}\int_{\Om} \mathbf{b} \times \mathbf{m}_{\epsilon} \cdot \vectorv  \dx \dt \\
 & \quad +\frac{1}{\epsilon^2} \int_0^{\tau}\int_{\Om} (\bar{\varrho}+ \epsilon q-\varrho_{\epsilon})  P^{\prime \prime}(\tilde{\varrho}) \; \partial_{t}(\bar{\varrho}+ \epsilon q) \dx\dt \\
 &\quad + \frac{1}{\epsilon^2} \int_0^{\tau}\int_{\Om}((\bar{\varrho}+ \epsilon q)\vectorv- \vectorm_{\epsilon}) \nabla_x P^{\prime  \prime}(\tilde{\varrho}) \nabla_{x}(\bar{\varrho}+ \epsilon q) \dx \dt\\
 & \quad - \int_0^{\tau}\int_{\Om} \nabla_x \tilde{\vectoru} : \text{d}[\mathfrak{R}_{v} +\mathfrak{R}_{p} \mathbb{I}]  (t,\cdot)\dt . \\
 \end{align*}

 Using the fact $\Divh\vectorv=0$, we obtain,
 \begin{align*}
 \Epsilon_{\epsilon}(\tau)+&  \int_{ \Om}\bigg[ \frac{1}{2} \text{d}\text{ trace} \mathfrak{R}_{v_{\epsilon}}(\tau,\cdot)+ \frac{1}{\gamma-1} \text{d}\mathfrak{R}_{p_{\epsilon}}(\tau,\cdot) \bigg]\\
 \quad &\leq \Epsilon_{\epsilon}(0)  -\int_0^{\tau}\int_{\Om}  (\vectorm_{\epsilon}-\varrho_{\epsilon} \vectorv )\cdot ( \partial_{t} \vectorv  + (\vectorv\cdot \nabla_{x_h})\vectorv )\dx \dt\\
 &\quad  - \int_0^{\tau}\int_{\Om}(\frac{ (\vectorm_{\epsilon}-\varrho_{\epsilon} \vectorv) \otimes( \vectorm_{\epsilon}-\varrho_{\epsilon}\vectorv)}{\varrho_{\epsilon}}) : \nabla_{x}\vectorv  \dx \dt  \\
 & \quad +\frac{1}{\epsilon} \int_0^{\tau}\int_{\Om} P^{\prime \prime}(\bar{\varrho}) \mathbf{m}_{\epsilon} \cdot \nabla_{x_h} q \dx \dt \\
 & \quad + \int_0^{\tau}\int_{\Om} ( q-\varrho_{\epsilon}^{1})  P^{\prime \prime}(\bar{\varrho}+ \epsilon q) \; \partial_{t}q \dx\dt \\
 &\quad - \frac{1}{\epsilon}\int_0^{\tau}\int_{\Om} \vectorm_{\epsilon}\cdot   P^{\prime  \prime}(\bar{\varrho}+ \epsilon q) \nabla_{x}q \dx \dt\\
 & \quad - \int_0^{\tau}\int_{\Om} \nabla_x \tilde{\vectoru} : \text{d}[\mathfrak{R}_{v} +\mathfrak{R}_{p} \mathbb{I}]  (t,\cdot)\dt =\Sigma_{i=1}^{7} \mathcal{L}_i .\\
 \end{align*}
  Consideration of \emph{well prepared data} yields,
  \begin{align*}
  \Epsilon_{\epsilon}(\varrho_{0,\epsilon}, \vectorm_{0,\epsilon} \;|\; \bar{\varrho}+ \epsilon q_0, \vectorv_0) \leq \Vert \frac{\vectorm_{0,\epsilon}}{\varrho_{0,\epsilon}}-  \vectorv_{0} \Vert_{L^2(\Om)}^2+ \Vert \varrho_{0,\epsilon}^{(1)}-q_0 \Vert^2_{L^2(\Om)}.
  \end{align*}
  Hence we conclude,
  \begin{align}\label{L1}
  \vert \mathcal{L}_1 \vert \leq c(\epsilon).
    \end{align}
    Here $c(\epsilon)$ is a generic function such that $c(\epsilon)\rightarrow 0$ as $\epsilon\rightarrow 0$.\\
  We write,
  $$\vectorm_{\epsilon}-\varrho_{\epsilon} \vectorv = (\vectorm_{\epsilon}-\bar{\varrho} \vectoru)+ ((\bar{\varrho}-\varrho_{\epsilon})\vectoru)+(\varrho_{\epsilon}(\vectoru-\vectorv))$$
  and $$ P^{\prime \prime}(\bar{\varrho}+ \epsilon q) =P^{\prime \prime}(\bar{\varrho}+ \epsilon q) -P^{\prime \prime}(\bar{\varrho}) +P^{\prime \prime}(\bar{\varrho}) . $$

  Now using above and results previous part, we obtain,
  \begin{align*}
  \vert \mathcal{L}_2 + \mathcal{L}_5 \vert \leq &c(\epsilon) + \int_0^{\tau}\int_{\Om}  \bar{\varrho}(\vectorv- \vectoru )\cdot ( \partial_{t} \vectorv  + (\vectorv\cdot \nabla_{x_h})\vectorv )\dx \dt \\
  & + \int_0^{\tau}\int_{\Om} ( q-\varrho^{(1)})  P^{\prime \prime}(\bar{\varrho}) \; \partial_{t}q \dx\dt.
  \end{align*}

  Now using properties of $q,\vectoru$ and $\vectorv$ we can conclude that,
  \begin{align*}
   \int_0^{\tau}\int_{\Om}  \bar{\varrho}(\vectorv- \vectoru )\cdot ( \partial_{t} \vectorv  + (\vectorv\cdot \nabla_{x_h})\vectorv )\dx \dt  + \int_0^{\tau}\int_{\Om} ( q-\varrho^{(1)})  P^{\prime \prime}(\bar{\varrho}) \; \partial_{t}q \dx\dt=0,
  \end{align*}
  since,
  \begin{align*}
  &\int_{\Om}  \bar{\varrho}(\vectorv- \vectoru )\cdot ( \partial_{t} \vectorv  + (\vectorv\cdot \nabla_{x_h})\vectorv )\dx \dt  + \int_0^{\tau}\int_{\Om} ( q-\varrho^{(1)})  P^{\prime \prime}(\bar{\varrho}) \; \partial_{t}q \dx\dt\\
  &=\frac{\text{d}}{\dt}\int_{ \Om}\big( \bar{\varrho}\vert \vectorv  \vert ^2 + P^{\prime\prime}(\bar{\varrho}) \vert q \vert^2 \big) \dx\\
  &-\int_{\Om}  \bar{\varrho} \vectoru \cdot ( \partial_{t} \vectorv  + (\vectorv\cdot \nabla_{x_h})\vectorv )\dx   - \int_{\Om} \varrho^{(1)}  P^{\prime \prime}(\bar{\varrho}) \; \partial_{t}q \dx\\
  &=\frac{p^{\prime}(\bar{\varrho})^2}{\bar{\varrho}}\frac{\text{d}}{\dt}\int_{ \Om}\big( \vert \nabla_{x_h} q  \vert ^2 + \frac{1}{p^{\prime}(\bar{\varrho})} \vert q \vert^2 \big) \dx\\
  &+ \frac{p^{\prime}(\bar{\varrho})^2}{\bar{\varrho}} \int_{\Om} \varrho^{(1)}\partial_{t}(\Delta_{x_h} q- \frac{1}{p^{\prime}(\bar{\varrho})} q)\dx- \bar{\varrho}\int_{ \Om} \vectoru \cdot (\vectorv\cdot \nabla_{x_h})\vectorv \dx\\
   &=\frac{p^{\prime}(\bar{\varrho})^2}{\bar{\varrho}}\frac{\text{d}}{\dt}\int_{ \Om}\big( \vert \nabla_{x_h} q  \vert ^2 + \frac{1}{p^{\prime}(\bar{\varrho})} \vert q \vert^2 \big) \dx\\
  &- \frac{p^{\prime}(\bar{\varrho})^2}{\bar{\varrho}} \int_{\Om} \varrho^{(1)}\nabla_{x_h}(\Delta_{x_h}q)\cdot \vectorv \dx- \bar{\varrho}\int_{ \Om} \vectoru \cdot (\vectorv\cdot \nabla_{x_h})\vectorv \dx\\
  &=\frac{p^{\prime}(\bar{\varrho})^2}{\bar{\varrho}}\frac{\text{d}}{\dt}\int_{ \Om}\big( \vert \nabla_{x_h} q  \vert ^2 + \frac{1}{p^{\prime}(\bar{\varrho})} \vert q \vert^2 \big) \dx\\
  &- {p^{\prime}(\bar{\varrho})} \int_{\Om} \varrho^{(1)}\nabla_{x_h}(\Delta_{x_h}q)\cdot \vectorv \dx- \bar{\varrho}\int_{ \Om} \vectoru \cdot (\vectorv\cdot \nabla_{x_h})\vectorv \dx\\
   &=\frac{p^{\prime}(\bar{\varrho})^2}{\bar{\varrho}}\frac{\text{d}}{\dt}\int_{ \Om}\big( \vert \nabla_{x_h} q  \vert ^2 + \frac{1}{p^{\prime}(\bar{\varrho})} \vert q \vert^2 \big) \dx\\
  &- {\bar{\varrho}} \int_{\Om} (\mathbf{b}\times \vectoru)(\Delta_{x_h}q)\cdot \vectorv \dx- \bar{\varrho}\int_{ \Om} \vectoru \cdot (\vectorv\cdot \nabla_{x_h})\vectorv \dx\\
  &=\frac{p^{\prime}(\bar{\varrho})^2}{\bar{\varrho}}\frac{\text{d}}{\dt}\int_{ \Om}\big( \vert \nabla_{x_h} q  \vert ^2 + \frac{1}{p^{\prime}(\bar{\varrho})} \vert q \vert^2 \big) \dx - \bar{\varrho} \int_{ \Om} \vectoru \cdot\nabla_{x_h} \vert \vectorv \vert^2 \dx.
  \end{align*}
  Clearly the third term vanishes. Now it is enough to prove,
  \begin{align*}
  \frac{\text{d}}{\dt}\int_{ \Om}\big( \vert \nabla_{x_h} q  \vert ^2 + \frac{1}{p^{\prime}(\bar{\varrho})} \vert q \vert^2 \big) \dx=0.
  \end{align*}
  The above inequality is true by multiplying second equation of target system by $q$.\\
  Hence we have,
  \begin{align}\label{L25}
   \vert \mathcal{L}_2 + \mathcal{L}_5 \vert \leq c(\epsilon).
  \end{align}
  Now we also obtain,
  \begin{align}\label{L3}
  \vert \mathcal{L}_{3} \vert \leq \int_{0}^{\tau} \Epsilon_{\epsilon} (t) \dt
  \end{align}
Next we consider,
  \begin{align*}
  \vert \mathcal{L}_4 + \mathcal{L}_{6} \vert & \leq \vert  \int_{0}^{\tau} \int_{ \Om} \frac{1}{\epsilon} \big( P^{\prime \prime }(\bar{\varrho} + \epsilon q) - P^{\prime \prime }(\bar{\varrho})\big) \vectorm_{\epsilon} \cdot (\nabla_{x_h} q, 0) \dx \dt \vert.
  \end{align*}
  It is easy to verify that,
   \begin{align*}
   \frac{1}{\epsilon} \big( P^{\prime \prime }(\bar{\varrho} + \epsilon q) - P^{\prime \prime }(\bar{\varrho})\big) \rightarrow P^{\prime \prime \prime }(\bar{\varrho})q \text{ in } L^{\infty}(0,T;L^{\infty} \cap L^2(\Om)).
   \end{align*}
  From here we can conclude that,
  \begin{align*}
  \mathcal{L}_4 + \mathcal{L}_{6} \rightarrow \int_{0}^{\tau} \int_{ \Om}  \bar{\varrho}  P^{\prime \prime \prime}(\bar{\varrho})q\vectorv \cdot (\nabla_{x_h}q,0) \dx \dt=0.
  \end{align*}
  Hence we have
  \begin{align}\label{L46}
  \vert \mathcal{L}_4 + \mathcal{L}_{6} \vert  \leq c(\epsilon).
  \end{align}
  Also we obtain,
  \begin{align}\label{L7}
  \vert \mathcal{L}_{7} \vert \leq  \int_0^{\tau}\int_{\Om} \mathbb{I} : \text{d}[\mathfrak{R}_{v} +\mathfrak{R}_{p} \mathbb{I}]  (t,\cdot)\dt .
  \end{align}

  Combining all these \eqref{L1}-\eqref{L7} we have,
  \begin{align}
  \begin{split}
  \Epsilon_{\epsilon}(\tau)+&  \int_{ \Om}\bigg[ \frac{1}{2} \text{d}\text{ trace} \mathfrak{R}_{v_{\epsilon}}(\tau,\cdot)+ \frac{1}{\gamma-1} \text{d}\mathfrak{R}_{p_{\epsilon}}(\tau,\cdot) \bigg]\\
  &\leq c(\epsilon)+  \int_{0}^{\tau} \Epsilon_{\epsilon} (t) \dt + C \int_{0}^{\tau} \int_{ \Om}\bigg[  \text{d}\text{ trace} \mathfrak{R}_{v_{\epsilon}}(\tau,\cdot)+ \text{d}\mathfrak{R}_{p_{\epsilon}}(\tau,\cdot) \bigg] \dt .
  \end{split}
  \end{align}
  Using Gronwall lemma we have,
  \begin{align}
  \begin{split}
  \Epsilon_{\epsilon}(\tau)+&  \int_{ \Om}\bigg[ \frac{1}{2} \text{d}\text{ trace} \mathfrak{R}_{v_{\epsilon}}(\tau,\cdot)+ \frac{1}{\gamma-1} \text{d}\mathfrak{R}_{p_{\epsilon}}(\tau,\cdot) \bigg] \\ &\leq c(\epsilon) C(T),
  \end{split}
  \end{align}
  where $c(\epsilon)\rightarrow 0$ as $\epsilon \rightarrow 0$.
  Hence,
  \begin{align}\label{lim_1_f}
  \lim\limits_{\epsilon\rightarrow 0} \Epsilon_{\epsilon} (\tau )=0.
  \end{align}
   Now using coerceivity of relative energy functional we can say,
   \begin{align}
   \varrho^{(1)}=q \text{ and } \vectorv=\vectoru.
   \end{align}
  From  \eqref{lim_1_f} we can conclude
  \begin{align*}
  &\frac{\varrho_{\epsilon}-\bar{\varrho}}{\epsilon} \rightarrow q \text{ \emph{strongly} in }L^1_{\text{loc}}((0,T)\times \Om),\\
  &\frac{\vectorm_{\epsilon}}{\sqrt{\varrho_{\epsilon}}} \rightarrow \sqrt{\bar{\varrho}}\vectorv \text{ \emph{strongly} in }L^1_{\text{loc}}((0,T)\times \Om;\R^3).
  \end{align*}
  It ends the proof of the theorem.

  \section{ Singular limit for "ill-prepared data"}
  \subsection{Definition of ill-prepared data}
  We define the \emph{ill-prepared initial data} as,
  \begin{Def}
  	We say that the set of initial data $(\varrho_{0,\epsilon},\vectorm_{0,\epsilon})_{(\epsilon>0)} $ is "well-prepared" if,
  	\begin{align}\label{ill_prepared_id}
  	\begin{split}
  	&\varrho_{0,\epsilon}=\bar{\varrho}+\epsilon \varrho_{0,\epsilon}^{(1)},\;\; \{\varrho_{0,\epsilon}\}_{\{\epsilon>0\} }\text{ is bounded in } L^2\cap L^{\infty}(\Om),\; \varrho_{0,\epsilon}^{(1)}\rightarrow \varrho_0^{(1)} \text{ in } L^2(\Om),\\
  	&\frac{\vectorm_{0,\epsilon}}{\varrho_{0,\epsilon}}\rightarrow  \vectorv_{0} \text{ in } L^2(\Omega;\R^3\big).
  	\end{split}
  	\end{align}
  \end{Def}

\subsection{Main Theorem:}
\begin{Th}
	Let pressure $p$ follows \eqref{p-condition}. Further we assume that the initial data is \emph{ill-prepared}, i.e. it follows \eqref{ill_prepared_id}. Also we assume intial data has better regularity,
	\begin{align}
	\varrho_{0}^{(1)} \in W^{k-1,2}(\Om),\; \vectoru_{0} \in W^{k,2}(\Om;\R^3),\; k\geq 3
	\end{align} Then after taking a subsequence, the following holds,
	\begin{align}
	\begin{split}
	&\varrho_{\epsilon}^{(1)} \equiv \frac{\varrho_{\epsilon}-\bar{\varrho}}{\epsilon} \rightarrow q \text{ weakly-(*) in } L^{\infty}(0,T;L^2+L^\gamma(\Om)),\\
	&\vectorm_\epsilon \rightarrow \bar{\varrho} \vectorv \text{ weakly-(*) in } L^{\infty}(0,T;L^2+L^{\frac{2\gamma}{\gamma+1}}(\Om)).\\
	\end{split}
	\end{align}
	Then $(q,\vectorv)$ solves \eqref{limit_system} with initial data $q(0,\cdot)=q_0$ where $q_0$ is the solution of elliptic equation,
	\begin{align}
	-\Delta_{x_h} q_0 + \frac{1}{p^{\prime}(\bar{\varrho})} q_0= \bar{\varrho} \int_{0}^{1} \curlh P_h(\vectoru_0) \; \text{d}x_3 + \frac{1}{p^{\prime}(\bar{\varrho})} \int_{0}^{1} \varrho_{0}^{1}\;																																																																																																																																																																	 \text{d}x_3.
	\end{align}
\end{Th}

One major problem for the analysis of the singular limit is the presence of rapidly oscillating Rossby-acoustic waves. Although these acoustic waves disappear in the course of low-mach number limit.
\subsubsection{Informal justification of derivation of the acoustic wave equation:} We rewrite the equation in fast time scale $\tau= \frac{t}{\epsilon}$ as,
\begin{align}
\frac{1}{\epsilon} \partial_{\tau} \varrho + \text{div}_x (\vectorm)&=0,
\end{align}
\begin{align}
\frac{1}{\epsilon} \partial_{\tau} \mathbf{m} + \Div (\frac{ \mathbf{m} \otimes \mathbf{m}}{\varrho})+\frac{1}{\epsilon^2}\nabla_x p(\varrho)+\frac{1}{\epsilon} \mathbf{b} \times \mathbf{m}&=0.
\end{align}
Again considering,
\begin{align*}
&\varrho_{\epsilon}= \bar{\varrho} + \epsilon \varrho_{\epsilon}^{(1)}+ \epsilon^2  \varrho_{\epsilon}^{(2)}+ \cdots,\\
&\vectorm_\epsilon= \bar{\varrho}\vectorv + \epsilon  \vectorm_\epsilon^{(1)} + \epsilon^2 \vectorm_\epsilon^{(2)} + \cdots.
\end{align*}
and comparing zeroth order term for continuity equation and $o(\epsilon^{-1})$ for momentum equation, we have,
\begin{align*}
\partial_{\tau} \varrho_{\epsilon}^{(1)} +\Div (\bar{\varrho} \vectorv) + o(\epsilon)=0,\\
\bar{\varrho} \partial_{\tau}(\vectorv)+ p^{\prime}(\bar{\varrho}) \nabla_{x}\varrho_{\epsilon}^{(1)}+ \bar{\varrho} \mathbf{b} \times \vectorv + o(\epsilon)=0.
\end{align*}

Going back to previous time scale we have,
\begin{align*}
\epsilon \partial_{t} r +\bar{\varrho} \Div  \vectorV =0,\\
\epsilon \partial_{t}\vectorV+ \frac{p^{\prime}(\bar{\varrho})}{\bar{\varrho}} \nabla_{x} r +  \mathbf{b} \times \vectorV =0.
\end{align*}
This motivates us how we obtain this hyperbolic system which represents acoustic waves.

\subsection{Proof of the theorem:}

The proof is in similar lines as in \emph{well-prepared case}
The first part is quite similar as we have done for the \emph{well-prepared data} in Convergence: part 1. Now for second part we need some extra term and dissipative estimate to deal it.
We recall the same relative energy ineqality,
\begin{align*}
\Epsilon_{\epsilon}(\tau)& +\int_{ \Om}\bigg[ \frac{1}{2} \text{d}\text{ trace} \mathfrak{R}_{v_{\epsilon}}(\tau,\cdot)+ \frac{1}{\gamma-1} \text{d}\mathfrak{R}_{p_{\epsilon}}(\tau,\cdot) \bigg]\\
\quad &\leq \Epsilon_{\epsilon}(0)  -\int_0^{\tau}\int_{\Om}  (\vectorm_{\epsilon}-\varrho_{\epsilon} \tilde{\vectoru})\cdot \partial_{t} \tilde{\vectoru}\dx \dt\\
&\quad  - \int_0^{\tau}\int_{\Om}(\frac{ (\vectorm_{\epsilon}-\varrho_{\epsilon} \tilde{\vectoru}) \otimes \vectorm_{\epsilon}}{\varrho_{\epsilon}}) : \nabla_x \tilde{\vectoru} \dx \dt  \\
&\quad  - \frac{1}{\epsilon^2} \int_0^{\tau}\int_{\Om} (p(\varrho_{\epsilon})-p(\tilde{\varrho})) \Div \tilde{\vectoru} \dx \dt \\
& \quad +\frac{1}{\epsilon} \int_0^{\tau}\int_{\Om} \mathbf{b} \times \mathbf{m}_{\epsilon} \cdot ( \tilde{\vectoru}- \frac{\vectorm_{\epsilon}}{\varrho_{\epsilon}})  \dx \dt  +\frac{1}{\epsilon^2} \int_0^{\tau}\int_{\Om} (\tilde{\varrho}-\varrho_{\epsilon}) \partial_{t} P^{\prime}(\tilde{\varrho}) \dx\dt \\
&\quad + \frac{1}{\epsilon^2} \int_0^{\tau}\int_{\Om}(\tilde{\varrho}\tilde{\vectoru}- \vectorm_{\epsilon}) \nabla_x P^{\prime}(\tilde{\varrho}) \dx \dt - \int_0^{\tau}\int_{\Om} \nabla_x \tilde{\vectoru} : \text{d}[\mathfrak{R}_{v} +\mathfrak{R}_{p} \mathbb{I}]  (t,\cdot)\dt . \\
\end{align*}
We rewrite the above one as,
\begin{align*}
\Epsilon_{\epsilon}(\tau)+& \int_{ \Om}\bigg[ \frac{1}{2} \text{d}\text{ trace} \mathfrak{R}_{v_{\epsilon}}(\tau,\cdot)+ \frac{1}{\gamma-1} \text{d}\mathfrak{R}_{p_{\epsilon}}(\tau,\cdot) \bigg] \\
\quad &\leq \Epsilon_{\epsilon}(0)  -\int_0^{\tau}\int_{\Om}  (\vectorm_{\epsilon}-\varrho_{\epsilon} \tilde{\vectoru})\cdot (\partial_{t} \tilde{\vectoru} + (\tilde{\vectoru}\cdot \nabla_{x}) \tilde{\vectoru})\dx \dt\\
&\quad  - \int_0^{\tau}\int_{\Om}(\frac{ (\vectorm_{\epsilon}-\varrho_{\epsilon} \tilde{\vectoru}) \otimes (\vectorm_{\epsilon}-\varrho_{\epsilon}\tilde{\vectoru})}{\varrho_{\epsilon}}) : \nabla_x \tilde{\vectoru} \dx \dt  \\
&\quad  - \frac{1}{\epsilon^2} \int_0^{\tau}\int_{\Om} (p(\varrho_{\epsilon})-p(\tilde{\varrho})) \Div \tilde{\vectoru} \dx \dt \\
& \quad +\frac{1}{\epsilon} \int_0^{\tau}\int_{\Om} \mathbf{b} \times \mathbf{m}_{\epsilon} \cdot ( \tilde{\vectoru}- \frac{\vectorm_{\epsilon}}{\varrho_{\epsilon}})  \dx \dt \\
& \quad +\frac{1}{\epsilon^2} \int_0^{\tau}\int_{\Om} (\tilde{\varrho}-\varrho_{\epsilon}) \partial_{t} P^{\prime}(\tilde{\varrho}) \dx\dt \\
&\quad + \frac{1}{\epsilon^2} \int_0^{\tau}\int_{\Om}(\tilde{\varrho}\tilde{\vectoru}- \vectorm_{\epsilon}) \nabla_x P^{\prime}(\tilde{\varrho}) \dx \dt\\
& \quad - \int_0^{\tau}\int_{\Om} \nabla_x \tilde{\vectoru} : \text{d}[\mathfrak{R}_{v} +\mathfrak{R}_{p} \mathbb{I}]  (t,\cdot)\dt . \\
\end{align*}
\subsubsection{ Initial data decomposition }
Following \cite{FN2014} subsection 4.3 we will regularize initial data $\varrho_{0}^{(1)},\; \vectoru_{0}=\frac{\vectorm_{0}}{\varrho_{0}}$ as, $$
\{[\varrho_{0}^{(1)}]_{\delta},\;
[\vectoru_{0}]_{\delta}\}_{\delta>0}, $$ with the help of cut-off function and Fourier transforms to gain both integrability and smoothness.
Now decompose
\begin{align*}
&[\varrho_{0}^{(1)}]_{\delta}=s_{0,\delta} + q_{0,\delta},\;[\vectoru_{0}]_{\delta}=\vectorV_{0,\delta}+ \vectorv_{0,\delta}
\end{align*}
with
\begin{align*}
&\text{ where, } -\Delta_{x_h} q_{0,\delta} + \frac{1}{p^{\prime}(\bar{\varrho})} q_{0,\delta}= \bar{\varrho} \int_{0}^{1} \curlh P_h[(\vectoru_0)]_{\delta} \; \text{d}x_3 + \frac{1}{p^{\prime}(\bar{\varrho})} \int_{0}^{1} [\varrho_{0}^{(1)}]_{\delta}\; \text{d}x_3\\
&\text{ and } \vectorv_{0,\delta}= \nabla_{x_h}^{\perp} q_{0,\delta}.
\end{align*}
Here we consider $s_{\epsilon,\delta},\; \vectorV_{\epsilon,\delta}$ solve
\begin{align*}
\epsilon \partial_{t} s_{\epsilon,\delta} +\bar{\varrho} \Div  \vectorV_{\epsilon,\delta} =0,\\
\epsilon \partial_{t}\vectorV_{\epsilon,\delta}+ \frac{p^{\prime}(\bar{\varrho})}{\bar{\varrho}} \nabla_{x} s_{\epsilon,\delta} +  \mathbf{b} \times \vectorV_{\epsilon,\delta} =0.
\end{align*}
with initial data $s_{0,\delta},\; \vectorV_{0,\delta} $.\\
Also cosider $q_\delta, \vectorv_{\delta}$ solve
\begin{align}
\begin{split}
&\frac{p^{\prime}(\bar{\varrho})}{\bar{\varrho}}\nabla_x q_{\delta} +  \mathbf{b} \times \vectorv_{\delta} =0,\\
&\partial_{t}(\Delta_{x_h}q_{\delta}-\frac{1}{p^{\prime}(\bar{\varrho})}q_{\delta}) + (\nabla_{x_h}^{\perp}q_{\delta}\cdot \nabla_{x_h})(\Delta_{x_h}q_{\delta}- \frac{1}{p^{\prime}(\bar{\varrho})}q_{\delta})=0,
\end{split}
\end{align}
with initial data
$q_{0,\delta},\vectorv_{0,\delta}$.

Next we recall the result from \cite{FN2014} Section 4 which is one of the important result concerning dispersive estimates of acoustis-rossby waves.
\begin{Prop}
	For the choice of initial data $s_{0,\delta},\; \vectorV_{0,\delta} $ as above we have the following result when  $\epsilon\rightarrow0$,
	\begin{align*}
	s_{\epsilon,\delta}\rightarrow 0 \text{ in } L^{p}(0,T;W^{m,\infty}(\Om)),\\
	\vectorV_{\epsilon,\delta} \rightarrow 0 \text{ in } L^{p}(0,T;W^{m,\infty}(\Om)),
	\end{align*}
	for any fixed  $\delta>0$, any $1 \leq p< \infty$ and $m=0,1,2,\cdots $.
\end{Prop}
\subsubsection{Relative Energy inequality and Convergence: Part 2}

First we have same apriori bounds and convergence of $\varrho_{\epsilon}, \vectorm_{\epsilon}$ as in \emph{well-prepared} case (See Convergence: Part 1). Now we need to choose a different test functions in this case.

Considering $$\tilde{\varrho}= \bar{\varrho}+ \epsilon (q_{\delta}+ s_{\epsilon,\delta}),\; \tilde{\vectoru} = \vectorv_{\delta}
+ \vectorV_{\epsilon,\delta},$$ we rewrite the previous inequality as,

\begin{align*}
\Epsilon_{\epsilon}(\tau)&+ \int_{ \Om}\bigg[ \frac{1}{2} \text{d}\text{ trace} \mathfrak{R}_{v_{\epsilon}}(\tau,\cdot)+ \frac{1}{\gamma-1} \text{d}\mathfrak{R}_{p_{\epsilon}}(\tau,\cdot) \bigg] \\
\quad &\leq \Epsilon_{\epsilon}(0)  -\int_0^{\tau}\int_{\Om}  (\vectorm_{\epsilon}-\varrho_{\epsilon} \tilde{\vectoru})\cdot (\partial_{t} \vectorv_{\delta} + (\tilde{\vectoru}\cdot \nabla_{x}) \tilde{\vectoru})\dx \dt\\
&\quad -\int_0^{\tau}\int_{\Om}  (\vectorm_{\epsilon}-\varrho_{\epsilon} \tilde{\vectoru})\cdot \partial_{t} \vectorV_{\epsilon,\delta} \dx \dt\\
&\quad  - \int_0^{\tau}\int_{\Om}(\frac{ (\vectorm_{\epsilon}-\varrho_{\epsilon} \tilde{\vectoru}) \otimes (\vectorm_{\epsilon}-\varrho_{\epsilon}\tilde{\vectoru})}{\varrho_{\epsilon}}) : \nabla_x \tilde{\vectoru} \dx \dt  \\
&\quad  - \frac{1}{\epsilon^2} \int_0^{\tau}\int_{\Om} (p(\varrho_{\epsilon})-p(\tilde{\varrho})-p^{\prime}(\tilde{\varrho})(\varrho_{\epsilon}-\tilde{\varrho})) \Div \tilde{\vectoru} \dx \dt \\
&\quad  + \frac{1}{\epsilon^2} \int_0^{\tau}\int_{\Om} (p^{\prime}(\tilde{\varrho})(\varrho_{\epsilon}-\tilde{\varrho})) \Div \tilde{\vectoru} \dx \dt \\
& \quad +\frac{1}{\epsilon} \int_0^{\tau}\int_{\Om} \mathbf{b} \times \mathbf{m}_{\epsilon} \cdot ( \tilde{\vectoru}- \frac{\vectorm_{\epsilon}}{\varrho_{\epsilon}})  \dx \dt +\frac{1}{\epsilon^2} \int_0^{\tau}\int_{\Om} (\tilde{\varrho}-\varrho_{\epsilon}) \partial_{t} P^{\prime}(\tilde{\varrho}) \dx\dt \\
&\quad + \frac{1}{\epsilon^2} \int_0^{\tau}\int_{\Om}(\tilde{\varrho}\tilde{\vectoru}- \vectorm_{\epsilon}) \nabla_x P^{\prime}(\tilde{\varrho}) \dx \dt- \int_0^{\tau}\int_{\Om} \nabla_x \tilde{\vectoru} : \text{d}[\mathfrak{R}_{v} +\mathfrak{R}_{p} \mathbb{I}]  (t,\cdot)\dt . \\
\end{align*}
Now,
\begin{align*}
&(p^{\prime}(\tilde{\varrho})(\varrho_{\epsilon}-\tilde{\varrho})) \Div \tilde{\vectoru}+ (\tilde{\varrho}-\varrho_{\epsilon}) \partial_{t} P^{\prime}(\tilde{\varrho})+(\tilde{\varrho}\tilde{\vectoru}- \vectorm_{\epsilon}) \nabla_x P^{\prime}(\tilde{\varrho})\\
=& \nabla_{x} P^{\prime}(\tilde{\varrho})(\varrho_{\epsilon} \tilde{\vectoru}-\vectorm_{\epsilon})+ \epsilon (\tilde{\varrho}-\varrho_{\epsilon}) P^{\prime \prime}(\tilde{\varrho}) (\partial_{t} q_{\delta}+ \Div[(q_{\delta}+s_{\epsilon,\delta})\tilde{\vectoru}])\\
=& \epsilon  P^{\prime\prime}(\tilde{\varrho})(\nabla_{x} q_{\delta}+s_{\epsilon,\delta})(\varrho_{\epsilon} \tilde{\vectoru}-\vectorm_{\epsilon})+ \epsilon (\tilde{\varrho}-\varrho_{\epsilon}) P^{\prime \prime}(\tilde{\varrho}) (\partial_{t} q_{\delta}+ \Div[(q_{\delta}+s_{\epsilon,\delta})\tilde{\vectoru}]).
\end{align*}
Further we have,
\begin{align*}
&\frac{1}{\epsilon} \mathbf{b} \times \vectorm_{\epsilon} \cdot \tilde{\vectoru} -(\vectorm_{\epsilon}-\varrho_{\epsilon} \tilde{\vectoru})\cdot \partial_{t} \vectorV_{\epsilon,\delta} + \frac{1}{\epsilon} P^{\prime\prime}(\tilde{\varrho})(\nabla_{x} q_{\delta}+s_{\epsilon,\delta})(\varrho_{\epsilon} \tilde{\vectoru}-\vectorm_{\epsilon})\\
=& \frac{1}{\epsilon}\big( P^{\prime\prime}(\bar{\varrho})-P^{\prime\prime}(\tilde{\varrho})\big) \nabla_{x} q_{\delta} \vectorm_{\epsilon} + \frac{1}{\epsilon}\big( P^{\prime\prime}(\bar{\varrho})-P^{\prime\prime}(\tilde{\varrho})\big) \nabla_{x} s_{\epsilon,\delta} \vectorm_{\epsilon} \\
& \quad - \frac{1}{\epsilon}\big( P^{\prime\prime}(\bar{\varrho})-P^{\prime\prime}(\tilde{\varrho})\big) \nabla_{x} s_{\epsilon,\delta} \varrho_{\; \epsilon} \tilde{\vectoru} - \frac{1}{\epsilon}\big( P^{\prime\prime}(\bar{\varrho})-P^{\prime\prime}(\tilde{\varrho})\big) \nabla_{x} q_{\delta} \varrho_{\epsilon} \tilde{\vectoru}\\
&= \frac{1}{\epsilon} \big( P^{\prime\prime}(\bar{\varrho})-P^{\prime\prime}(\tilde{\varrho})\big) \nabla_{x}(q_\delta + s_{\epsilon,\delta}) (\vectorm_{\epsilon}-\varrho_{\epsilon} \tilde{\vectoru}).
	\end{align*}
	
Hence we have,
\begin{align*}
\Epsilon_{\epsilon}(\tau)+&\int_{ \Om}\bigg[ \frac{1}{2} \text{d}\text{ trace} \mathfrak{R}_{v_{\epsilon}}(\tau,\cdot)+ \frac{1}{\gamma-1} \text{d}\mathfrak{R}_{p_{\epsilon}}(\tau,\cdot) \bigg]\\
\quad &\leq \Epsilon_{\epsilon}(0)  -\int_0^{\tau}\int_{\Om}  (\vectorm_{\epsilon}-\varrho_{\epsilon} \tilde{\vectoru})\cdot (\partial_{t} \vectorv_{\delta} + (\tilde{\vectoru}\cdot \nabla_{x}) \tilde{\vectoru})\dx \dt\\
&\quad  - \int_0^{\tau}\int_{\Om}(\frac{ (\vectorm_{\epsilon}-\varrho_{\epsilon} \tilde{\vectoru}) \otimes (\vectorm_{\epsilon}-\varrho_{\epsilon}\tilde{\vectoru})}{\varrho_{\epsilon}}) : \nabla_x \tilde{\vectoru} \dx \dt  \\
&\quad  - \frac{1}{\epsilon^2} \int_0^{\tau}\int_{\Om} (p(\varrho_{\epsilon})-p(\tilde{\varrho})-p^{\prime}(\tilde{\varrho})(\varrho_{\epsilon}-\tilde{\varrho})) \Div \tilde{\vectoru} \dx \dt \\
&\quad + \int_{0}^{\tau} \int_{ \Om} \frac{\tilde{\varrho}- \varrho_{\epsilon}}{\epsilon}P^{\prime \prime}(\tilde{\varrho}) (\partial_{t} q_{\delta}+ \Div[(q_{\delta}+s_{\epsilon,\delta})\tilde{\vectoru}]) \dx \dt \\
& \quad + \int_{0}^{\tau} \int_{ \Om} \frac{1}{\epsilon}\big( P^{\prime\prime}(\bar{\varrho})-P^{\prime\prime}(\tilde{\varrho})\big) \nabla_{x}(q_\delta + s_{\epsilon,\delta}) (\vectorm_{\epsilon}-\varrho_{\epsilon} \tilde{\vectoru})\dx \dt\\
& \quad- \int_0^{\tau}\int_{\Om} \nabla_x \tilde{\vectoru} : \text{d}[\mathfrak{R}_{v} +\mathfrak{R}_{p} \mathbb{I}]  (t,\cdot)\dt  \\
&\quad = \sum_{i=1}^{7} \mathcal{T}_i
\end{align*}
We introduce a notaion $h_i(\epsilon,\delta)$ as $\lim\limits_{\epsilon \rightarrow 0} h_i(\epsilon,\delta)= \bar{h}_i(\delta)$ and $ \lim\limits_{\delta\rightarrow 0} \bar{h}_i(\delta)=0$ for $i\in \mathbb{N}$.
\begin{align*}
\vert \mathcal{T}_1 \vert \leq h_1(\epsilon,\delta).
\end{align*}
\begin{align*}
 \mathcal{T}_2 = & \int_0^{\tau}\int_{\Om}  (\vectorm_{\epsilon}-\varrho_{\epsilon} \tilde{\vectoru})\cdot (\partial_{t} \vectorv_{\delta} + (\tilde{\vectoru}\cdot \nabla_{x}) \tilde{\vectoru})\dx \dt\\
 & =\int_0^{\tau}\int_{\Om}  (\vectorm_{\epsilon}-\varrho_{\epsilon} \vectorv_{\delta})\cdot (\partial_{t} \vectorv_{\delta} + (\vectorv_{\delta}\cdot \nabla_{x}) \vectorv_{\delta})\dx \dt + \bar{\mathcal{T}}_2.
\end{align*}
\begin{align*}
\vert \bar{\mathcal{T}}_2 \vert \leq h_2(\epsilon,\delta), \text{ with } \bar{h}_2(\delta)=0.
\end{align*}
Eventually we have,
\begin{align*}
&\int_0^{\tau}\int_{\Om}  (\vectorm_{\epsilon}-\varrho_{\epsilon} \vectorv_{\delta})\cdot (\partial_{t} \vectorv_{\delta} + (\vectorv_{\delta}\cdot \nabla_{x}) \vectorv_{\delta})\dx \dt \\
 & \rightarrow \int_0^{\tau}\int_{\Om}  \bar{\varrho}(\vectoru- \vectorv_{\delta})\cdot (\partial_{t} \vectorv_{\delta} + (\vectorv_{\delta}\cdot \nabla_{x}) \vectorv_{\delta})\dx \dt
\end{align*}
\begin{align*}
\vert \mathcal{T}_3 \vert + \vert \mathcal{T}_4 \vert \leq  \sup_{  t\in [0,T]} \Vert \vectorv_\delta(t,\cdot ) \Vert_{W^{1,\infty}(\mathbb{R}^2;\mathbb{R}^2)} \int_{0}^{\tau} \Epsilon_{\epsilon}(t) \dt.
\end{align*}

For $\mathcal{T}_5$, we use similar treatment as $\mathcal{T}_2$ we have
\begin{align*}
\mathcal{T}_5=\int_{0}^{\tau} \int_{ \Om} \frac{\tilde{\varrho}- \varrho_{\epsilon}}{\epsilon}P^{\prime \prime}({\tilde{\varrho}}) (\partial_{t} q_{\delta}+ \Div[q_\delta \vectorv_{\delta}]) \dx \dt +  \bar{\mathcal{T}_5},
\end{align*}
with $\vert \bar{\mathcal{T}}_5 \vert \leq h_5(\epsilon,\delta), \text{ with } \bar{h}_5(\delta)=0$ and
\begin{align*}
&\int_{0}^{\tau} \int_{ \Om} \frac{\tilde{\varrho}- \varrho_{\epsilon}}{\epsilon}P^{\prime \prime}({\tilde{\varrho}}) (\partial_{t} q_{\delta}+ \Div[q_\delta \vectorv_{\delta}]) \dx \dt \\
&\rightarrow \int_{0}^{\tau} \int_{ \Om}(q_\delta - \varrho^{(1)})P^{\prime \prime}({\bar{\varrho}})(\partial_{t} q_{\delta}+ \Div[q_\delta \vectorv_{\delta}]) \dx \dt.
\end{align*}
 Following similar reasoning line as in \emph{well-prepared} case, we have,
 \begin{align}
 \begin{split}
& \int_0^{\tau}\int_{\Om}  \bar{\varrho}(\vectoru- \vectorv_{\delta})\cdot (\partial_{t} \vectorv_{\delta} (\vectorv_{\delta}\cdot \nabla_{x}) \vectorv_{\delta})\dx \dt \\
& \quad  + \int_{0}^{\tau} \int_{ \Om}(q_\delta - \varrho^{(1)})(\partial_{t} q_{\delta}+ \Div[q_\delta \vectorv_{\delta}]) \dx \dt=0.
 \end{split}
 \end{align}
 for $\mathcal{T}_6$  we have,
 \begin{align*}
 \mathcal{T}_6=\int_{0}^{\tau} \int_{ \Om} \frac{1}{\epsilon}\big( P^{\prime\prime}(\bar{\varrho})-P^{\prime\prime}(\tilde{\varrho})\big) \nabla_{x}(q_\delta ) (\vectorm_{\epsilon}-\varrho_{\epsilon} \vectorv_{\delta})\dx \dt + \bar{\mathcal{T}_6},
 \end{align*}
 with $\vert \bar{\mathcal{T}}_6 \vert \leq h_6(\epsilon,\delta), \text{ with } \bar{h}_6(\delta)=0,$ and
 \begin{align*}
 &\int_{0}^{\tau} \int_{ \Om} \frac{1}{\epsilon}\big( P^{\prime\prime}(\bar{\varrho})-P^{\prime\prime}(\tilde{\varrho})\big) \nabla_{x}(q_\delta ) (\vectorm_{\epsilon}-\varrho_{\epsilon} \vectorv_{\delta})\dx \dt\\
 &\rightarrow \int_{0}^{\tau} \int_{ \Om} P^{\prime \prime \prime }(\bar{\varrho})
q_\delta \nabla_{x}q_\delta \bar{\varrho} (\vectoru-\vectorv_{\delta}) \dx \dt=0.
 \end{align*}
We also have,
\begin{align*}
\vert \mathcal{T}_7 \vert  \leq \sup_{  t\in [0,T]} \Vert \vectorv_\delta(t,\cdot ) \Vert_{W^{1,\infty}(\mathbb{R}^2;\mathbb{R}^2)}\int_0^{\tau }\int_{ \Om}\bigg[ \frac{1}{2} \text{d}\text{ trace} \mathfrak{R}_{v_{\epsilon}}(\tau,\cdot)+ \frac{1}{\gamma-1} \text{d}\mathfrak{R}_{p_{\epsilon}}(t,\cdot) \bigg]
\end{align*}
Thus we have
\begin{align*}
\Epsilon_{\epsilon}&(\tau)+\int_{ \Om}\bigg[ \frac{1}{2} \text{d}\text{ trace} \mathfrak{R}_{v_{\epsilon}}(\tau,\cdot)+ \frac{1}{\gamma-1} \text{d}\mathfrak{R}_{p_{\epsilon}}(\tau,\cdot) \bigg]\\
& \quad \leq h(\epsilon,\delta) + C(q) \bigg( \int_{0}^{\tau} \Epsilon_{\epsilon}(t) \dt \\
& \quad \quad \quad \quad \quad \;+ \int_0^{\tau }\int_{ \Om}\bigg[ \frac{1}{2} \text{d}\text{ trace} \mathfrak{R}_{v_{\epsilon}}(\tau,\cdot)+ \frac{1}{\gamma-1} \text{d}\mathfrak{R}_{p_{\epsilon}}(t,\cdot) \bigg]\bigg)
\end{align*}
Passing limit for first for $\epsilon\rightarrow 0$, then for $\delta \rightarrow 0$ we can conclude the proof of theorem.

\section{Concluding Remarks:}
Clearly we can see \emph{well-prepared} case as a particuler case of \emph{ill prepared} case. Although here we show that for well prepared case you dont need any dispersive estimates.\\

\vspace{5mm}

\centerline{ \bf Acknowledgement}
\vspace{2mm}

The work is supported by Einstein Stiftung, Berlin. I would like to thank my  Ph.D supervisor Prof. E. Feireisl for his valuable suggestions and comments.



\end{document}